# Near-optimal analysis of Lasserre's univariate measure-based bounds for multivariate polynomial optimization

Lucas Slot[1] · Monique Laurent[1,2]

**Abstract**
We consider a hierarchy of upper approximations for the minimization of a polynomial $f$ over a compact set $K \subseteq \mathbb{R}^n$ proposed recently by Lasserre (arXiv:1907.097784, 2019). This hierarchy relies on using the push-forward measure of the Lebesgue measure on $K$ by the polynomial $f$ and involves univariate sums of squares of polynomials with growing degrees $2r$. Hence it is weaker, but cheaper to compute, than an earlier hierarchy by Lasserre (SIAM Journal on Optimization 21(3), 864–885, 2011), which uses multivariate sums of squares. We show that this new hierarchy converges to the global minimum of $f$ at a rate in $O(\log^2 r / r^2)$ whenever $K$ satisfies a mild geometric condition, which holds, eg., for convex bodies and for compact semialgebraic sets with dense interior. As an application this rate of convergence also applies to the stronger hierarchy based on multivariate sums of squares, which improves and extends earlier convergence results to a wider class of compact sets. Furthermore, we show that our analysis is near-optimal by proving a lower bound on the convergence rate in $\Omega(1/r^2)$ for a class of polynomials on $K = [-1, 1]$, obtained by exploiting a connection to orthogonal polynomials.

This work is supported by the European Union's EU Framework Programme for Research and Innovation Horizon 2020 under the Marie Skłodowska-Curie Actions Grant Agreement No 764759 (MINOA).

✉ Monique Laurent
  Monique.Laurent@cwi.nl

  Lucas Slot
  Lucas.Slot@cwi.nl

[1] Centrum Wiskunde & Informatica (CWI), Amsterdam, The Netherlands

[2] Tilburg University, Tilburg, The Netherlands

Published online: 30 October 2020

## 1 Introduction

Consider the problem of finding the minimum value taken by an $n$-variate polynomial $f \in \mathbb{R}[x]$ over a compact set $K \subseteq \mathbb{R}^n$, i.e., computing the parameter:

$$f_{\min} = \min_{x \in K} f(x). \quad (1)$$

Throughout we also set $f_{\max} = \max_{x \in K} f(x)$. Computing the parameter $f_{\min}$ (or $f_{\max}$) is a hard problem in general, including for instance the maximum stable set problem as a special case. For a general reference on polynomial optimization and its applications, we refer, eg., to [14,16].

If we fix a Borel measure $\lambda$ with support $K$, problem (1) may be reformulated as minimizing the integral $\int_K f(x)\sigma(x)d\lambda(x)$ over all sum-of-squares polynomials $\sigma \in \Sigma[x]$ that provide a probability density on $K$ with respect to the measure $\lambda$. By bounding the degree of $\sigma$, we obtain the following hierarchy of upper bounds on $f_{\min}$ proposed by Lasserre [15]:

$$f_{\min} \leq f^{(r)} := \min \left\{ \int_K f(x)\sigma(x)d\lambda(x) : \sigma \in \Sigma[x]_r, \int_K \sigma(x)d\lambda(x) = 1 \right\}. \quad (2)$$

Here $\Sigma[x]$ denotes the set of polynomials that can be written as a sum of squares of polynomials and we set $\Sigma[x]_r = \Sigma[x] \cap \mathbb{R}[x]_{2r}$. Since sums of squares of polynomials can be expressed using semidefinite programming, for any fixed $r \in \mathbb{N}$ the parameter $f^{(r)}$ can be computed efficiently by semidefinite programming or, even simpler, as the smallest eigenvalue of an appropriate matrix of size $\binom{n+r}{r}$ ([15], see also [5]).

Recently, Lasserre [17] introduced new, weaker but more economical, upper bounds on $f_{\min}$ that are based on a *univariate* approach to the problem. For this purpose, he considers the *push-forward measure* $\lambda_f$ of $\lambda$ by $f$, which is defined by

$$\lambda_f(B) = \lambda(f^{-1}(B)) \quad \text{for any Borel set } B \subseteq \mathbb{R}. \quad (3)$$

Note that for any measurable function $g : \mathbb{R} \to \mathbb{R}$, we thus have

$$\int_{f(K)} g(t)d\lambda_f(t) = \int_K g(f(x))d\lambda(x). \quad (4)$$

We then can define the following hierarchy of upper bounds on $f_{\min}$:

$$\begin{aligned} f_{\min} \leq f_{\text{pfm}}^{(r)} &:= \min \left\{ \int_{f(K)} ts(t)d\lambda_f(t) : s \in \Sigma[t]_r, \int_{f(K)} s(t)d\lambda_f(t) = 1 \right\} \\ &= \min \left\{ \int_K f(x)s(f(x))d\lambda(x) : s \in \Sigma[t]_r, \int_K s(f(x))d\lambda(x) = 1 \right\}. \end{aligned} \quad (5)$$

The difference with the parameter $f^{(r)}$ is that we now restrict the search to *univariate* sums of squares $s \in \Sigma[t]_r$, which we then evaluate at the polynomial $f$, leading to the multivariate sum of squares $\sigma_{\text{pfm}} := s \circ f \in \Sigma[x]_{rd}$ if $f$ has degree $d$. Therefore we have the inequality

$$f_{\min} \leq f^{(rd)} \leq f^{(r)}_{\text{pfm}}. \tag{6}$$

Again, the parameter $f^{(r)}_{\text{pfm}}$ can be computed efficiently for any fixed $r$. But now it can be computed as the smallest eigenvalue of an appropriate matrix of much smaller size $r+1$ (see (8) below). Asymptotic convergence of the parameters $f^{(r)}_{\text{pfm}}$ to $f_{\min}$ is shown in [17], but no quantitative results are given there. In this paper, we are interested in analyzing the convergence rate of the parameters $f^{(r)}_{\text{pfm}}$ to the global minimum $f_{\min}$ in terms of the degree $r$.

### 1.1 Previous work

In what follows we always consider for $\lambda$ the Lebesgue measure on $K$ (unless specified otherwise). Several results exist on the convergence rate of the parameters $f^{(r)}$ to the global minimum $f_{\min}$, depending on the set $K$. The best rates in $O(1/r^2)$ were shown in [5,6,23] when $K$ belongs to special classes of convex bodies, including the hypercube $[-1, 1]^n$, the ball $B^n$, the sphere $S^{n-1}$, the standard simplex $\Delta^n$ and compact sets that are locally 'ball-like'. Furthermore, it was shown in [5] that this analysis is best possible in general (already for $K = [-1, 1]$ and $f(x) = x$). The starting point for each of these results is a connection between the parameters $f^{(r)}$ and the smallest roots of certain orthogonal polynomials (see [5, Sect. 2] and the short recap below).

In [23, Theorems 10-11], a rate in $O(\log^2 r/r^2)$ was shown for general convex bodies $K$, as well as a rate in $O(\log r/r)$ for general compact sets $K$ that satisfy a minor geometric condition (a srengthening of Assumption 1 below). There the analysis relied on constructing explicit sum-of-squares densities that approximate well the Dirac delta function at a global minimizer of $f$, making use of the so-called 'needle' polynomials from [12]. An improved rate in $O(\log^k r/r^k)$ was shown in [23, Theorem 14] when the partial derivatives of $f$ up to degree $k-1$ vanish at one of its global minimizers on $K$.

When $K$ is a convex body, a convergence rate in $O(1/r)$ had been shown earlier in [4], by exploiting a link to simulated annealing. There the authors considered sum-of-squares densities of (roughly) the form $\sigma = s \circ f$, where $s(t) = \sum_{k=0}^{2r}(-t/T)^k/k! \in \Sigma[t]_r$ is the truncated Taylor expansion of the exponential $e^{-t/T}$. Hence this specific choice of $s$ (or $\sigma$) provides an upper bound not only for the parameter $f^{(rd)}$ (as exploited in [4]) but also for the parameter $f^{(r)}_{\text{pfm}}$ and thus the result of [4] gives directly $f^{(r)}_{\text{pfm}} - f_{\min} = O(1/r)$ when $K$ is a convex body.

The result above gives a first quantitative analysis of the parameters $f^{(r)}_{\text{pfm}}$ for convex bodies. In this paper we improve this result in two directions. First we sharpen the analysis and show the stronger convergence rate $O(\log^2 r/r^2)$ and second we show that this analysis applies a large class of compact sets (those satisfying Assumption 1), which includes all semialgebraic sets that have a dense interior.

We also mention briefly another hierarchy of bounds when $K$ is a semialgebraic set, of the form

$$K = \{x \in \mathbb{R}^n : g_1(x) \geq 0, \ldots, g_m(x) \geq 0\} \text{ with } g_1, \ldots, g_m \in \mathbb{R}[x].$$

Then lower bounds for the minimum of $f$ over $K$ can be obtained as

$$f_{(r)} = \sup \left\{ \lambda : f - \lambda = \sum_{j=0}^{m} \sigma_j g_j, \right.$$

$$\left. \sigma_j \text{ sum-of-squares polynomials, } \deg(\sigma_j g_j) \leq 2r \ (0 \leq j \leq m) \right\}$$

(setting $g_0 = 1$). This hierarchy has been widely studied in the literature (see, eg., [14,16] and references therein). Asymptotic convergence to $f_{\min}$ holds when the semi-algebraic set $K$ satisfies the Archimedean condition (which implies $K$ is compact) [13] and relies on the positivity certificate of Putinar [20]. (The Archimedean condition requires existence of $R > 0$ such that $R - \sum_{i=1}^{n} x_i^2$ lies in the quadratic module generated by the $g_j$'s, consisting of the polynomials $\sum_j s_j g_j$ for some sum-of-squares polynomials $s_j$). The question arises naturally of analyzing the quality of the bounds $f_{(r)}$. A convergence rate in $O(1/(\log(r/c))^{1/c})$ is shown in [18], where $c$ is a constant depending only on $K$. If in the definition of the bounds $f_{(r)}$ we allow decompositions in the preordering, which consists of the polynomials $\sum_{J \subseteq [m]} \sigma_J \prod_{j \in J} g_j$ with $\sigma_J$ sum-of-squares polynomials, then, based on Schmüdgen's positivity certificate [21], asymptotic convergence holds for any compact $K$ and a stronger convergence rate in $O(1/r^c)$ was shown in [22] (where $c$ again depends only on $K$). When allowing decompositions in the preordering a stronger convergence rate in $O(1/r)$ was shown for special sets like the simplex (in [1]) and the hypercube (in [2]). (See also [7] for an overview). For the minimization of a homogeneous polynomial $f$ over the unit sphere an improved convergence rate in $O(1/r^2)$ for the bounds $f_{(r)}$ was shown recently in [11] (improving the earlier rate in $O(1/r)$ from [9]). It turns out that this analysis relies (implicitly) on the convergence rate of the upper bounds for a special class of polynomials. This indicates there are intimate links between the upper and lower bounds $f^{(r)}$ and $f_{(r)}$, which forms an additional motivation for better understanding the upper bounds $f^{(r)}$. Showing an improved convergence analysis for the bounds $f_{(r)}$ for broader classes of semialgebraic sets remains an important research question.

### 1.2 New results

The main contribution of this paper is the following bound on the convergence rate of the parameter $f_{\text{pfm}}^{(r)}$ that holds whenever $K$ satisfies a minor geometric condition.

**Theorem 1** *Let $K \subseteq \mathbb{R}^n$ be a compact connected set satisfying Assumption 1 below. Then we have*

$$f^{(r)}_{\text{pfm}} - f_{\min} = O(\log^2 r / r^2).$$

In view of (6), we immediately get the following corollary, extending the rate in $O(\log^2 r/r^2)$, shown in [23] for convex bodies, to all connected compact sets $K$ satisfying Assumption 1.

**Corollary 1** *Let $K \subseteq \mathbb{R}^n$ be a compact connected set satisfying Assumption 1. Then we have*

$$f^{(r)} - f_{\min} = O(\log^2 r / r^2).$$

In light of the following special case of [5, Corollary 3.2] our result on the convergence rate of $f^{(r)}_{\text{pfm}}$ is best possible in general, up to the log-factor.

**Theorem 2** [5] *Let $K = [-1, 1]$ and let $f(x) = x$. Then $f^{(r)} = -1 + \Theta(1/r^2)$. As a direct consequence, we have $f^{(r)}_{\text{pfm}} (= f^{(r)}) = -1 + \Omega(1/r^2)$.*

As an additional result, we extend the lower bound $\Omega(1/r^2)$ on the error range $f^{(r)}_{\text{pfm}} - f_{\min}$ to the class of functions $f(x) = x^{2k}$ with integer $k \geq 1$.

**Theorem 3** *Let $K = [-1, 1]$ and let $f(x) = x^{2k}$ for $k \geq 1$ integer. Then we have $f^{(r)}_{\text{pfm}} = \Omega(1/r^2)$.*

Combining Theorem 3 with the fact that $f^{(r)} = O(\log^{2k} r / r^{2k})$ when $f(x) = x^{2k}$ (using [23, Theorem 14]), we thus show a large separation between the asymptotic quality of the bounds $f^{(r)}$ and $f^{(r)}_{\text{pfm}}$ for this class of functions.

### 1.3 Approach and discussion

As already mentioned above, a crucial ingredient in the analysis of the parameters $f^{(r)}$ for special compact sets like the hypercube $[-1, 1]^n$, the ball, the sphere, or the simplex, is the analysis in the univariate case when $K = [-1, 1]$ (equipped with the Lebesgue measure or more generally allowing a weight of Jacobi type) and the special polynomial $f(t) = t$. Let $\{p_i \in \mathbb{R}[t]_i : i \in \mathbb{N}\}$ be the (unique) orthonormal basis of $\mathbb{R}[t]$ with respect to the inner product $\langle \cdot, \cdot \rangle_\lambda$ given by

$$\langle p, q \rangle_\lambda = \int_K pq \, d\lambda \quad \text{for } p, q \in \mathbb{R}[t]. \tag{7}$$

Then, as is shown in [5], the parameter $f^{(r)}$ coincides with the smallest eigenvalue of the (truncated) *moment matrix* $M_{\lambda,r}$ of $\lambda$, which is defined as

$$M_{\lambda,r} := \left( \int_K t p_i p_j \, d\lambda \right)^r_{ij=0}. \tag{8}$$

A classical result on orthogonal polynomials (cf., eg., [24]) shows that the eigenvalues of $M_{\lambda,r}$ are given by the roots of $p_{r+1}$. Hence, the parameter $f^{(r)}$ is equal to the smallest root of $p_{r+1}$, the asymptotic behaviour of which is well understood and known to be in $-1 + \Theta(1/r^2)$ when $\lambda$ is a measure of Jacobi type ([5], see also Lemma 2 below).

Recall that $\lambda_f$ is the push-forward measure of $\lambda$ by $f$, as defined in (3), and $f(K) = [f_{\min}, f_{\max}]$ since we assume $K$ is compact and connected. Let $\{p_{f,i} : i \in \mathbb{N}\}$ denote the orthonormal basis of $\mathbb{R}[t]$ with respect to the inner product $\langle \cdot, \cdot \rangle_{\lambda_f}$ on the interval $[f_{\min}, f_{\max}]$. In view of the above discussion, if we use the first (univariate) formulation of $f_{\text{pfm}}^{(r)}$ in (5), we can immediately conclude that $f_{\text{pfm}}^{(r)}$ is equal to the smallest eigenvalue of the matrix

$$M_{\lambda_f,r} := \left( \int_{f_{\min}}^{f_{\max}} t p_{f,i} p_{f,j} d\lambda_f \right)_{i,j=0}^{r},$$

and also to the smallest root of the orthogonal polynomial $p_{f,r+1}$. However it is not clear how to exploit this connection in order to gain information about the convergence rate of the parameters $f_{\text{pfm}}^{(r)}$ since the orthogonal polynomials $p_{f,i}$ are not known explicitly in general.

In this paper, we will go back to the idea of trying to find a good sum-of-squares polynomial approximation of the Dirac delta function. As in [23], we make use of the needle polynomials from [12] for this purpose. The difference with the approach in [23] is that we now work on the interval $[f_{\min}, f_{\max}]$; so we need an approximation of the Dirac delta function centered at $f_{\min}$, which is on the boundary of this interval. As is already noted in [12], this special setting allows for better approximations than would be available in general.

### 1.4 Outline

The rest of the paper is organized as follows. In Sect. 2 we give a proof of Theorem 1. Then, in Sect. 3, we prove Theorem 3. We provide some numerical examples that illustrate the practical behaviour of the bounds $f^{(r)}$ and $f_{\text{pfm}}^{(r)}$ in Sect. 4. Finally, in Sect. 5, we give a small discussion of the geometric Assumption 1 below and we show that it is satisfied by the compact semialgebraic sets with a dense interior.

## 2 Convergence analysis for the new hierarchy

We first state the precise geometric condition alluded to in Theorem 1.

**Assumption 1** There exist positive constants $\epsilon_K, \eta_K > 0$ and $N \geq n$, such that, for all $x \in K$ and $0 < \delta \leq \epsilon_K$, we have

$$\text{vol}\left(K \cap B_\delta^n(x)\right) \geq \eta_K \delta^N \text{vol}\left(B^n\right). \tag{9}$$

Here, $B_\delta^n(x)$ is the Euclidean ball centered at $x$ with radius $\delta$ and $B^n = B_1^n(0)$.

A slightly stronger version of Assumption 1 (requiring $N = n$) was introduced in [3], where it was used to give the first error analysis in $O(1/\sqrt{r})$ for the bounds $f^{(r)}$. The condition of [3] is satisfied, eg., when $K$ is a convex body, or more generally when $K$ satisfies an interior cone condition, or when $K$ is star-shaped with respect to a ball (see also [3] for a more complete discussion). The weaker condition (9) is satisfied additionally by the compact semialgebraic sets that have a dense interior, which allows in particular that $K$ has certain types of cusps. We discuss Assumption 1 in more detail in Sect. 5 below.

We show the following restatement of Theorem 1.

**Theorem 4** *Assume $K$ is connected, compact and satisfies the above geometric condition (9). Then there exists a constant $C$ (depending only on $n$, the Lipschitz constant of $f$ and $K$) such that*

$$f^{(r)}_{\text{pfm}} - f_{\min} \leq C \frac{\log^2 r}{r^2} (f_{\max} - f_{\min}) \quad \text{for all large } r.$$

The rest of this section is devoted to the proof of Theorem 4. We will make the following assumptions in order to simplify notation in our arguments. Let $a$ be a global minimizer of $f$ in $K$. After applying a suitable translation (replacing $K$ by $K - a$ and the polynomial $f$ by the polynomial $x \mapsto f(x - a)$), we may assume that $a = 0$, that is, we may assume that the global minimum of $f$ over $K$ is attained at the origin. Furthermore, it suffices to work with the rescaled polynomial

$$F(x) := \frac{f(x) - f_{\min}}{f_{\max} - f_{\min}},$$

which satisfies $F(K) = [0, 1]$, with $F_{\min} = 0$ and $F_{\max} = 1$. Indeed, one can easily check that

$$f^{(r)}_{\text{pfm}} - f_{\min} \leq (f_{\max} - f_{\min}) F^{(r)}_{\text{pfm}}.$$

Then, for this polynomial $F$, we know that the support of the push-forward measure $\lambda_F$ is equal to $[0, 1]$, and (5) gives

$$\begin{aligned} F^{(r)}_{\text{pfm}} &= \min \left\{ \int_0^1 ts(t) d\lambda_F(t) : s \in \Sigma[t]_r, \int_0^1 s(t) d\lambda_F(t) = 1 \right\} \\ &= \min \left\{ \int_K F(x) s(F(x)) d\lambda(x) : s \in \Sigma[t]_r, \int_K s(F(x)) d\lambda(x) = 1 \right\}. \end{aligned} \quad (10)$$

In order to analyze the bound $F^{(r)}_{\text{pfm}}$, we follow a similar strategy to the one employed in [23] to analyze the bound $F^{(r)}$. Namely, we construct a univariate sum-of-squares polynomial $s$ which approximates well the Dirac delta centered at the origin on the interval $[0, 1]$, making use of the so-called $\frac{1}{2}$-*needle polynomials* from [12].

**Lemma 1** [12] *Let $h \in (0, 1)$ be a scalar and let $r \in \mathbb{N}$. Then there exists a univariate polynomial $v_r^h \in \Sigma[t]_{2r}$ satisfying the following properties:*

$$\begin{aligned} v_r^h(0) &= 1, \\ 0 \leq v_r^h(t) &\leq 1 \quad \text{for all } t \in [0, 1], \\ v_r^h(t) &\leq 4e^{-\frac{1}{2}r\sqrt{h}} \quad \text{for all } t \in [h, 1]. \end{aligned} \quad (11)$$

We consider the sum-of-squares polynomial $s(t) := Cv_r^h(t)$, where $h \in (0, 1)$ will be chosen later, and $C$ is chosen so that $s$ is a density on $[0, 1]$ with respect to the measure $\lambda_F$. That is,

$$C = \left( \int_K v_r^h(F(x)) d\lambda(x) \right)^{-1}.$$

As $s$ is a feasible solution to (10), we obtain

$$F_{\text{pfm}}^{(r)} \leq \int_K F(x) s(F(x)) d\lambda(x) = \frac{\int_K F(x) v_r^h(F(x)) d\lambda(x)}{\int_K v_r^h(F(x)) d\lambda(x)}.$$

Our goal is thus to show that

$$\text{ratio} := \frac{\int_K F(x) v_r^h(F(x)) d\lambda(x)}{\int_K v_r^h(F(x)) d\lambda(x)} = O\left( \frac{\log^2 r}{r^2} \right). \quad (12)$$

Define the set

$$K_h = \{x \in K : F(x) \leq h\}.$$

We first work out the numerator of (12), which we split into two terms, depending whether we integrate on $K_h$ or on its complement:

$$\begin{aligned} \int_K F(x) v_r^h(F(x)) d\lambda(x) &= \int_{K_h} F(x) v_r^h(F(x)) d\lambda(x) + \int_{K \setminus K_h} F(x) v_r^h(F(x)) d\lambda(x) \\ &\leq h \int_{K_h} v_r^h(F(x)) d\lambda(x) + \int_{K \setminus K_h} v_r^h(F(x)) d\lambda(x). \end{aligned}$$

Here we have upper bounded $F(x)$ by $h$ on $K_h$ and by 1 on $K \setminus K_h$. On the other hand, we can lower bound the denominator in (12) as follows:

$$\int_K v_r^h(F(x)) d\lambda(x) \geq \int_{K_h} v_r^h(F(x)) d\lambda(x).$$

Combining the above two inequalities on numerator and denominator we get

$$\text{ratio} \le h + \frac{\int_{K \setminus K_h} v_r^h(F(x)) d\lambda(x)}{\int_{K_h} v_r^h(F(x)) d\lambda(x)}.$$

Thus we only need to upper bound the second term above. We first work on the numerator. For any $x \in K \setminus K_h$ we have $F(x) > h$ and thus, using (11), we get $v_r^h(F(x)) \le 4e^{-\frac{1}{2}r\sqrt{h}}$. This implies

$$\int_{K \setminus K_h} v_r^h(F(x)) d\lambda(x) \le 4e^{-\frac{1}{2}r\sqrt{h}} \lambda(K).$$

Next, we bound the denominator. In [23, Corollary 4], it is observed that

$$v_r^h(t) \ge 1 - 32r^2 t \ge \frac{1}{2} \quad \text{for all } t \in [0, \frac{1}{64r^2}].$$

Set $\rho = \frac{1}{64r^2}$. We will later choose $h \ge \rho$, so that $K_h \supseteq K_\rho := \{x \in K : F(x) \le \rho\}$ and $v_r^h(F(x)) \ge \frac{1}{2}$ for all $x \in K_\rho$. As $K$ is compact, there exists a Lipschitz constant $C_F > 0$ such that

$$F(x) \le C_F \|x\| \quad \text{for all } x \in K. \tag{13}$$

Note that $K \cap B^n_{\rho/C_F} \subseteq K_\rho$. By the geometric assumption (9) we have

$$\lambda(K \cap B^n_{\rho/C_F}) \ge \eta_K \left(\frac{\rho}{C_F}\right)^N \lambda(B^n)$$

for all $r$ large enough such that $\rho/C_F \le \epsilon_K$. We can then lower bound the denominator as follows:

$$\int_{K_h} v_r^h(F(x)) d\lambda(x) \ge \int_{K_\rho} v_r^h(F(x)) d\lambda(x) \ge \frac{1}{2}\lambda(K_\rho) \ge \frac{1}{2}\lambda(K \cap B^n_{\rho/C_F})$$

$$\ge \frac{1}{2}\eta_K \left(\frac{\rho}{C_F}\right)^N \lambda(B^n).$$

Combining the above inequalities, we obtain

$$\text{ratio} \le h + \frac{e^{-\frac{1}{2}r\sqrt{h}}}{\rho^N} \cdot \frac{8 \cdot \lambda(K)(C_F)^N}{\eta_K \lambda(B^n)}.$$

If we now select $h = \left(4(N+1)\frac{\log r}{r}\right)^2$, we have $h \ge \rho$ and a straightforward computation shows that

$$\text{ratio} \le O\left(\frac{\log^2 r}{r^2}\right).$$

Here, the constant in the big O depends on $n$, $N$, $C_F$, $\eta_K$ and $\lambda(K)$. This concludes the proof of Theorem 4.

## 3 Separation for a special class of polynomials

In this section we consider in more detail the behaviour of the bounds $f^{(r)}$ and $f_{\text{pfm}}^{(r)}$ for the class of polynomials $f(x) = x^{2k}$ (with $k \geq 1$ integer) on the interval $K = [-1, 1]$. Then $f([-1, 1]) = [0, 1]$ and, by applying (6) to the polynomial $f(x) = x^{2k}$, we have the following inequality:

$$0 \leq f^{(2rk)} \leq f_{\text{pfm}}^{(r)} \quad \text{for any } r \geq 1.$$

Note that for any $i \leq 2k - 1$, the $i$th derivative of $f$ vanishes at its global minimizer 0 on $[-1, 1]$. Using [23, Theorem 14], we therefore have that $f^{(2rk)} \leq f^{(r)} = O(\log^{2k} r / r^{2k})$. On the other hand, the convergence rate in $O(\log^2 r / r^2)$ for $f_{\text{pfm}}^{(r)}$ shown in Theorem 1 is optimal up to the log-factor. Indeed, we will show here a lower bound for $f_{\text{pfm}}^{(r)}$ in $\Omega(1/r^2)$.

Let $\lambda_k := \lambda_f$ denote the push-forward measure (3) of the Lebesgue measure on $[-1, 1]$ by the function $f(x) = x^{2k}$, and let $\{p_{k,i}(t) : i \in \mathbb{N}\} \subseteq \mathbb{R}[t]$ denote the family of orthogonal polynomials that provide an orthonormal basis for $\mathbb{R}[t]$ w.r.t. the inner product $\langle \cdot, \cdot \rangle_{\lambda_k}$ (cf. (7)). Then, as shown in [5] and as recalled above, the parameter $f_{\text{pfm}}^{(r)}$ is equal to the smallest root of the polynomial $p_{k,r+1}(t)$. As it turns out, here we can find explicitly the push-forward measure $\lambda_k$, which can be shown to be of Jacobi type. Hence, we have information about the corresponding orthogonal polynomials $p_{k,i}$, whose extremal roots are well understood. First we introduce the classical Jacobi polynomials (see, eg., [24] for a general reference).

**Lemma 2** *Let $a, b > -1$. Consider the weight function $w_{a,b}(x) = (1 - x)^a (1 + x)^b$ on the interval $[-1, 1]$ and let $\{p_i^{a,b}(x) : i \in \mathbb{N}\}$ be the corresponding family of orthogonal polynomials. Then $p_i^{a,b}$ is known as the degree $i$ Jacobi polynomial (with parameters $a, b$), and its smallest root $\xi_i^{a,b}$ satisfies:*

$$\xi_i^{a,b} = -1 + \Theta(1/i^2). \tag{14}$$

**Proof** A proof of this fact based on results in [8,10] is given in [5]. □

**Lemma 3** *For any integrable function $g$ on $[-1, 1]$ we have the identity*

$$\int_{-1}^{1} g(x^{2k}) dx = \frac{1}{k} \int_{0}^{1} g(t) t^{-1+1/2k} dt.$$

*Hence, the push-forward measure $\lambda_k$ is given by $d\lambda_k(t) := \frac{1}{k} t^{-1+\frac{1}{2k}} dt$ for $t \in [0, 1]$.*

**Table 1** Polynomial test functions

| Name | Formula | $f_{\min}$ |
|---|---|---|
| Booth | $f_{bo}(x) = (10x_1 + 20x_2 - 7)^2 + (20x_1 + 10x_2 - 5)^2$ | $f_{bo}(\frac{1}{10}, \frac{3}{10}) = 0$ |
| Matyas | $f_{ma}(x) = 26(x_1^2 + x_2^2) - 48x_1 x_2$ | $f_{ma}(0, 0) = 0$ |
| Camel | $f_{ca}(x) = 50x_1^2 - \frac{2625}{4}x_1^4 + \frac{15625}{6}x_1^6 + 25x_1 x_2 + 25x_2^2$ | $f_{ca}(0, 0) = 0$ |
| Motzkin | $f_{mo}(x) = 64x_1^4 x_2^2 + 64x_1^2 x_2^4 - 48x_1^2 x_2^2 + 1$ | $f_{mo}(\pm\frac{1}{2}, \pm\frac{1}{2}) = 0$ |

In each case, $f_{\min}$ is the global minimum of $f$ on $[-1, 1]^2$

**Proof** It suffices to show the first claim, which follows by making a change of variables $t = x^{2k}$ so that we get

$$\int_{-1}^{1} g\left(x^{2k}\right) dx = 2 \int_{0}^{1} g\left(x^{2k}\right) dx = 2 \int_{0}^{1} g(t) \frac{t^{-1+\frac{1}{2k}}}{2k} dt = \frac{1}{k} \int_{0}^{1} g(t) t^{-1+\frac{1}{2k}} dt.$$

□

**Proof of Theorem 3** By applying the change of variables $x = 2t - 1$, we see that the Jacobi type measure $(1 - x)^a (1 + x)^b dx$ on $[-1, 1]$ corresponds to the measure $2^{a+b}(1 - t)^a t^b dt$ on $[0, 1]$ and that (up to scaling) the orthogonal polynomials for the latter measure on $[0, 1]$ are given by $t \mapsto p_i^{a,b}(2t - 1)$ for $i \in \mathbb{N}$.

If we set $a = 0$ and $b = -1 + 1/2k$, then the measure obtained in this way on $[0, 1]$ is precisely the push-forward measure $\lambda_k$ (see Lemma 3). Hence, we can conclude that (up to scaling) the orthogonal polynomials $p_{k,i}$ for $\lambda_k$ on $[0, 1]$ are given by $p_{k,i}(t) = p_i^{a,b}(2t - 1)$ for each $i \in \mathbb{N}$. Therefore, the smallest root of $p_{k,r+1}(t)$ is equal to $(\xi_{r+1}^{a,b} + 1)/2 = \Theta(1/r^2)$ by (14). In particular, we can conclude that $f_{\text{pfm}}^{(r)} = \Omega(1/r^2)$ for any $k \geq 1$.

□

## 4 Numerical examples

In this section, we illustrate the practical behaviour of the bounds $f_{\text{pfm}}^{(r)}$ and $f^{(r)}$ using some numerical examples. Recall from (6) that $f^{(dr)} \leq f_{\text{pfm}}^{(r)}$ if $f$ has degree $d$. However, we will see that in many of the examples below the parameter $f_{\text{pfm}}^{(r)}$ provides in fact a better bound than $f^{(r)}$.

**Comparison of $f_{\text{pfm}}^{(r)}$ and $f^{(r)}$ for polynomial test functions.** First, we consider the polynomial test functions listed in Table 1. These are all well-known in optimization, and were already used to test the behaviour of the bounds $f^{(r)}$ in [3,23]. We compare the bounds $f^{(r)}$ and $f_{\text{pfm}}^{(r)}$ directly for $f \in \{f_{bo}, f_{ma}, f_{ca}, f_{mo}\}$, computed for the unit box $[-1, 1]^2$ and the unit ball $B^2$. For $1 \leq r \leq 20$, we compute the values of the fraction:

$$\rho_r(f) := \frac{f_{\text{pfm}}^{(r)} - f_{\min}}{f^{(r)} - f_{\min}}.$$

So, values of $\rho_r(f)$ smaller than 1 indicate good performance of the bounds $f_{\text{pfm}}^{(r)}$ in comparison to $f^{(r)}$. The results can be found in Fig. 3. Remarkably, it appears that the performance of the bound $f_{\text{pfm}}^{(r)}$ is comparable to (or better than) the performance of $f^{(r)}$ in each instance, except for the Camel function. Additionally, we note that the performance of $f_{\text{pfm}}^{(r)}$ for the Motzkin polynomial is comparatively much better on the unit ball than on the unit box. Figure 1 shows a plot of the Camel function, as well as the sum-of-squares densities corresponding to $f^{(6)}$ and $f_{\text{pfm}}^{(6)}$ on the unit box. Note that while the density corresponding to $f^{(6)}$ resembles the Dirac delta function centered at the global minimizer (0, 0) of the Camel function, the density corresponding to $f_{\text{pfm}}^{(6)}$ instead mirrors the Camel function itself.

**Comparison of $f_{\text{pfm}}^{(r)}$ and $f^{(r)}$ for the special class of polynomials $f(x) = x^{2k}$.** Next, we consider the polynomials $f(x) = x^{2k}$ for $k \geq 1$ on the interval $[-1, 1]$, which were treated in Sect. 3. In Fig. 4, the values of $\rho_r(f)$ are shown for $1 \leq r \leq 20$ and $1 \leq k \leq 5$. It can be seen that the performance of $f_{\text{pfm}}^{(r)}$ is comparable to the performance of $f^{(r)}$ for $k = 1$ (indeed, in this case we have $f_{\text{pfm}}^{(r)} = f^{(2r)}$), but it is much worse for $k > 1$, which matches our earlier findings (Theorem 3). In Fig. 2, the optimal sum-of-squares densities $\sigma$ (corresponding to $f^{(r)}$) and $\sigma_{\text{pfm}}$ (corresponding to $f_{\text{pfm}}^{(r)}$) are depicted for $k = 1, 3, 5$ and $r = 6$. Note that while the density $\sigma$ changes very little as we increase $k$, the density $\sigma_{\text{pfm}}$ grows increasingly 'flat' around the minimizer 0 of $f$ (mirroring the behavior of $f$ itself). As such, the density $\sigma_{\text{pfm}}$ is a comparatively much worse approximation of the Dirac delta function centered at 0 than $\sigma$. Note also that in this instance $f^{(r)} = f^{(r+1)}$ for even $r$, explaining the 'zig-zagging' behaviour of the ratio $\rho_r(f)$.

**Comparison of $f_{\text{pfm}}^{(r)}$ and $f^{(r)}$ for random instances of maximum cut.** Finally, we consider some polynomial maximization problems on $[-1, 1]^n$ coming from small instances of MAXCUT. An instance of MAXCUT with vertex set $[n]$ and edge weights $w_{ij} \geq 0$ can be written as:

$$\text{OPT} := \max_{x \in [-1,1]^n} f(x), \text{ where } f(x) := \frac{1}{4} \sum_{i,j \in [n]} w_{ij}(x_i - x_j)^2. \quad (15)$$

Note that while $f$ is usually maximized over the discrete cube $\{-1, 1\}^n$, the formulation (15) is equivalent as $f$ is convex.

Following [15], we create our instances by setting $w_{ij} = 0$ with probability $p$, and sampling $w_{ij}$ uniformly from [0, 1] otherwise. In Table 2, we list values of $f_{\text{pfm}}^{(r)}$ and $f^{(r)}$ for a few such random instances with $p = 1/2$ and $n = 8$. In each case, $f_{\text{pfm}}^{(r)}$ provides a better bound than $f^{(r)}$. In Table 3, we list the average over 50 randomly generated instances of the ratios:

$$\text{ratio} = \frac{\text{OPT} - f^{(r)}}{\text{OPT}} \quad \text{and} \quad \text{ratio}_{\text{pfm}} = \frac{\text{OPT} - f_{\text{pfm}}^{(r)}}{\text{OPT}}$$

**Table 2** Values of $f^{(r)}$ and $f^{(r)}_{\text{pfm}}$ for randomly generated instances of MAXCUT ($n = 8$, $p = 1/2$)

|  | $r = 1$ | $r = 2$ | $r = 3$ | $r = 4$ | $r = 5$ | $r = 6$ | OPT |
|---|---|---|---|---|---|---|---|
| Ex1 |  |  |  |  |  |  |  |
| $f^{(r)}$ | 1.58 | 2.06 | 2.45 | 2.82 | 3.16 | 3.46 | 6.09 |
| $f^{(r)}_{\text{pfm}}$ | 1.98 | 2.65 | 3.20 | 3.67 | 4.06 | 4.38 |  |
| Ex2 |  |  |  |  |  |  |  |
| $f^{(r)}$ | 2.20 | 2.77 | 3.27 | 3.73 | 4.16 | 4.54 | 8.07 |
| $f^{(r)}_{\text{pfm}}$ | 2.61 | 3.41 | 4.11 | 4.72 | 5.25 | 5.71 |  |
| Ex3 |  |  |  |  |  |  |  |
| $f^{(r)}$ | 2.03 | 2.56 | 3.02 | 3.43 | 3.81 | 4.14 | 7.24 |
| $f^{(r)}_{\text{pfm}}$ | 2.46 | 3.19 | 3.81 | 4.34 | 4.79 | 5.15 |  |
| Ex4 |  |  |  |  |  |  |  |
| $f^{(r)}$ | 1.59 | 2.05 | 2.44 | 2.81 | 3.13 | 3.42 | 5.80 |
| $f^{(r)}_{\text{pfm}}$ | 1.98 | 2.62 | 3.15 | 3.60 | 3.97 | 4.28 |  |

**Table 3** Average performance of the bounds $f^{(r)}$ and $f^{(r)}_{\text{pfm}}$ for random instances of MAXCUT ($n = 8$)

|  | $r = 1$ | $r = 2$ | $r = 3$ | $r = 4$ |
|---|---|---|---|---|
| $p = 1/4$ |  |  |  |  |
| Ratio | 0.72 | 0.65 | 0.59 | 0.53 |
| Ratio$_{\text{pfm}}$ | 0.66 | 0.57 | 0.48 | 0.40 |
| $p = 1/2$ |  |  |  |  |
| Ratio | 0.73 | 0.65 | 0.59 | 0.53 |
| Ratio$_{\text{pfm}}$ | 0.68 | 0.56 | 0.47 | 0.39 |
| $p = 3/4$ |  |  |  |  |
| Ratio | 0.73 | 0.64 | 0.57 | 0.49 |
| Ratio$_{\text{pfm}}$ | 0.65 | 0.53 | 0.43 | 0.35 |

for $r \leq 4$ and $p \in \{1/4, 1/2, 3/4\}$. Although it seems $f^{(r)}_{\text{pfm}}$ is more sensitive to changes in the density of the instances, we find again that it provides a better bound in general than $f^{(r)}$.

## 5 On the geometric assumption

As mentioned above, the condition (9) is a weaker version of a condition introduced in [3]. There, the authors demand that there exist constants $\eta_K, \epsilon_K$ such that

$$\text{vol}\left(K \cap B^n_\delta(x)\right) \geq \eta_K \delta^n \, \text{vol}(B^n) \quad \forall x \in K, \ \forall 0 < \delta \leq \epsilon_K. \tag{16}$$

The difference is that the power of $\delta$ in (16) is fixed to be the dimension $n$ of $K$, whereas it is allowed to be an arbitrary $N \geq n$ in (9).

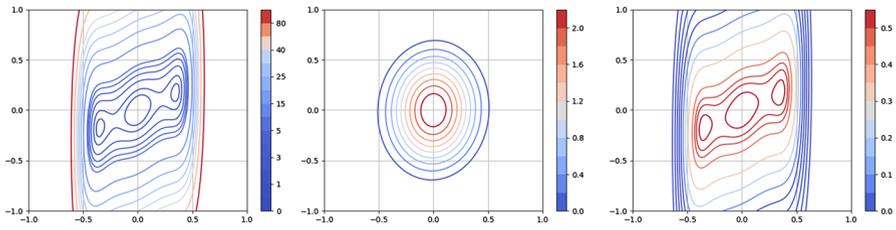

**Fig. 1** The Camel function (left) and its sum-of-squares densities corresponding to $f^{(6)}$ (middle) and $f_{\text{pfm}}^{(6)}$ (right) on the unit box

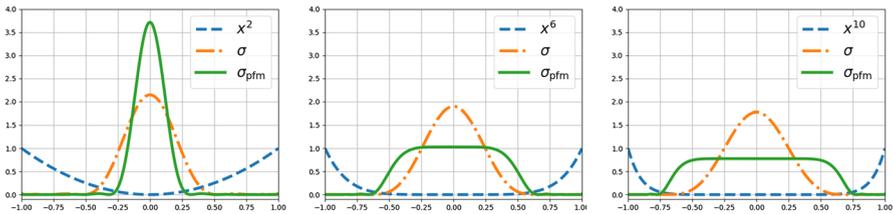

**Fig. 2** The functions $f(x) = x^{2k}$ and their sum-of-squares densities corresponding to $f^{(6)}$ and $f_{\text{pfm}}^{(6)}$ on the interval $[-1, 1]$ for $k = 1$ (left), $k = 3$ (middle) and $k = 5$ (right)

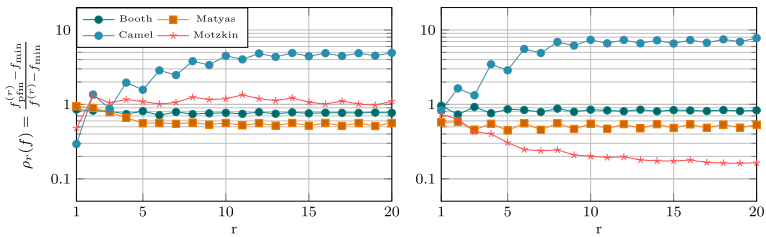

**Fig. 3** Comparison of the bounds $f^{(r)}$ and $f_{\text{pfm}}^{(r)}$ for the first four functions in Table 1, computed on the unit box (left) and unit ball (right)

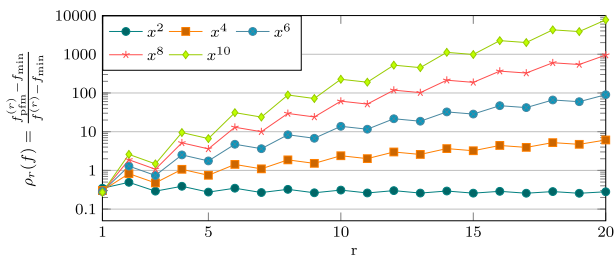

**Fig. 4** Comparison of the bounds $f^{(r)}$ and $f_{\text{pfm}}^{(r)}$ for functions of the form $f(x) = x^{2k}$ on the interval $[-1, 1]$

Condition (9) is satisfied by a significantly larger class of sets $K$ than (16). In particular, as we will observe below, sets satisfying (9) may have *polynomial* cusps, whereas sets satisfying (16) may not have any cusps at all.

**Example 1** Consider the set $K = \{x \in \mathbb{R}^2 : 0 \leq x_1 \leq 1,\ 0 \leq x_2 \leq x_1^2\}$ (see Fig. 5). This set $K$ satisfies (9) (with $N = 3$), but it does not satisfy (16). Indeed, for the point $0 \in K$ we have

$$\text{vol}\left(K \cap B_\delta^2(0)\right) \leq \int_0^\delta t^2 dt = \frac{1}{3}\delta^3,$$

and conclude (16) cannot be satisfied at $x = 0$. Note that the point 0 is indeed a polynomial cusp of the set $K$.

**Example 2** Consider the set $K = \{x \in \mathbb{R}^2 : 0 \leq x_1 \leq 1,\ 0 \leq x_2 \leq \exp(-1/x_1)\}$ (see Fig. 5). This set $K$ does not satisfy (9) (and, as a consequence, does not satisfy (16)). Indeed, for the point $0 \in K$ we have

$$\text{vol}\left(K \cap B_\delta^2(0)\right) \leq \int_0^\delta \exp(-1/t) dt \leq \delta \exp(-1/\delta).$$

Now note that for any $N, \eta > 0$ fixed, we have:

$$\lim_{\delta \to 0} \frac{\eta \delta^N}{\delta \exp(-1/\delta)} = \infty,$$

and so (9) can not be satisfied at $x = 0$. Note that the point 0 is an exponential cusp of $K$.

It turns out that compact semialgebraic sets which have a dense interior (aka being fat) satisfy Assumption 1, as is shown essentially in [19].

**Definition 1** A set $K \subseteq \mathbb{R}^n$ is called *fat* if $K \subseteq \overline{\text{int }K}$, i.e., the interior of $K$ is dense in $K$.

**Theorem 5** ([19], Theorem 6.4, see also Remark 6.5) *Let $K \subseteq \mathbb{R}^n$ be a compact, fat semialgebraic[1] set. Then there exist constants $\eta > 0$, $N \geq 1$ and a positive integer $d \in \mathbb{N}$ such that one may find a polynomial $h_x$ of degree $d$ for each $x \in K$ satisfying:*

$$h_x(0) = x, \tag{17}$$
$$h_x(t) \in K \text{ for } t \in [0, 1], \text{ and} \tag{18}$$
$$B_{\eta t^N}^n(h_x(t)) \subseteq K \text{ for } t \in [0, 1]. \tag{19}$$

*Furthermore, the polynomials $h_x$ may be chosen such that $\|x - h_x(t)\| \leq t$ for all $x \in K, t \in [0, 1]$.*

---

[1] In fact, the result is shown for *subanalytic* sets, of which semialgebraic sets are an example.

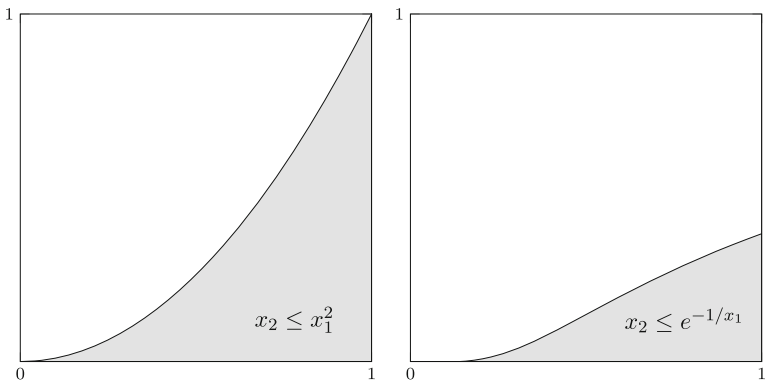

**Fig. 5** The set $K = \{x \in \mathbb{R}^2 : 0 \leq x_1 \leq 1, 0 \leq x_2 \leq x_1^2\}$ (left) and $K = \{x \in \mathbb{R}^2 : 0 \leq x_1 \leq 1, 0 \leq x_2 \leq \exp(-1/x_1)\}$ (right)

**Corollary 2** *Let $K \subseteq \mathbb{R}^n$ be a compact, fat semialgebraic set. Then $K$ satisfies Assumption 1.*

*Proof* For $x \in K$, let $\eta$, $N$ and $h_x$ be as in Theorem 5. We may assume that $h_t := h_x(t) \in B_t^n(x)$ for all $t \in [0, 1]$. For clarity, we write $B(y, a) := B_a^n(y)$ in the rest of the proof.

Using the triangle inequality and (19) we find that

$$\operatorname{vol} B\left(h_t, \eta t^N\right) \leq \operatorname{vol}\left(B\left(x, t + \eta t^N\right) \cap K\right) \leq \operatorname{vol}\left(B(x, (1 + \eta)t) \cap K\right)$$

for all $t \in [0, 1]$, noting that $t^N \leq t$ in this case. But now substituting $\delta = (1 + \eta)t$ yields:

$$\operatorname{vol}(B(x, \delta) \cap K) \geq \operatorname{vol} B\left(h_t, \eta \delta^N (1 + \eta)^{-N}\right) = \left(\frac{\eta}{(1+\eta)^N}\right)^n \cdot \delta^{Nn} \operatorname{vol} B(x, 1),$$

showing (9) holds for $0 < \delta \leq \epsilon_K := (1 + \eta)$ and $\eta_K = \eta^n (1 + \eta)^{-Nn}$. □

## 6 Conclusions

We have shown a convergence rate in $O(\log^2 r / r^2)$ for the approximations $f_{\text{pfm}}^{(r)}$ of the minimum of a polynomial $f$ over a compact connected set $K$ satisfying the minor geometric assumption (9). Furthermore, we have shown that this analysis is near-optimal, in the sense that the asymptotic behaviour of the error range $f_{\text{pfm}}^{(r)} - f_{\min}$ is in $O(\log^2 r / r^2)$ in general and in $\Omega(1/r^2)$ for an infinite class of polynomials.

This latter result shows that although the worst-case guarantees on the convergence of the bounds $f^{(r)}$ and $f_{\text{pfm}}^{(r)}$ are very similar, a large separation may exist for certain polynomials (eg., when $f(x) = x^{2k}$). Of course, it should be noted that the parameter $f_{\text{pfm}}^{(r)}$ can be obtained via a much smaller eigenvalue computation than the parameter $f^{(r)}$, namely by computing the smallest eigenvalue of a matrix of size $r + 1$ for the latter in comparison to a matrix of size $\binom{n+r}{r}$ for the former.

From a computational point of view, one should also observe that while the computation of $f_{\text{pfm}}^{(r)}$ involves a smaller matrix, it however requires to know the moments $\int_K f^k d\lambda$ of powers of $f$ for $k \leq 2r$. If $f$ has many terms this computation can be demanding for large values of $r$. This has to be taken into account when comparing the computational burden of both bounds $f^{(r)}$ and $f_{\text{pfm}}^{(r)}$.

Lastly, as a surprising consequence of Theorem 1, we are able to extend the bound in $O(\log^2 r/r^2)$ on the convergence rate of $f^{(r)}$ to all compact connected sets $K$ satisfying the geometric condition (9), whereas it was previously only known for convex bodies [23]. In this sense, the arguments of Sect. 2 can be seen as a refinement (and simplification) of the ones given in [23].

As said above, the analysis in this paper is near-optimal: we can show an upper bound in $O(\log^2 r/r^2)$ and a lower bound in $\Omega(1/r^2)$ for a certain class of polynomials. Deciding what is the right regime and whether the log-factor can be avoided in the convergence analysis is the main research question left open by this work.

The log-factor arises from our analysis technique, based on using polynomial approximation by the needle polynomials. We had to use this analysis technique since the behaviour of the orthogonal polynomials for the push-forward measure $\lambda_f$ is not known for general $f$. On the other hand, our results may be interpreted as giving back some information for general push-forward measures $\lambda_f$ and their corresponding orthogonal polynomials $p_{f,i}$ on the interval $[f_{\min}, f_{\max}]$. Indeed, what our results imply is that for any polynomial $f$ and any compact connected $K$ satisfying (9), the asymptotic behaviour of the smallest root of $p_{f,i}$ is in $f_{\min} + O(\log^2 r/r^2)$.

**Acknowledgements** We wish to thank Edouard Pauwels for bringing to our attention the needle polynomials, as well as the anonymous referees for their helpful suggestions.